\documentclass{article}

\usepackage{amsmath, amssymb}

\title{A Generalized Grover/Zeta Correspondence} 

\author{
Takashi KOMATSU \\
Math. Research Institute Calc for Industry, \\
Minami, Hiroshima, 732-0816, JAPAN \\ 
e-mail: ta.komatsu@sunmath-calc.co.jp \\ 
\\ 
Norio KONNO \\
Department of Applied Mathematics, Faculty of Engineering, \\ 
Yokohama National University \\
Hodogaya, Yokohama 240-8501, JAPAN \\
e-mail: konno-norio-bt@ynu.ac.jp \\ 
\\ 
Iwao SATO \\ 
Oyama National College of Technology,  \\
Oyama, Tochigi 323-0806, JAPAN \\ 
e-mail: isato@oyama-ct.ac.jp \\ 
\\ 
Shunya TAMURA \\
Graduate School of Science and Engineering,  \\
Yokohama National University,  \\ 
Hodogaya, Yokohama, 240-8501, JAPAN \\ 
e-mail: tamura-shunya-kj@ynu.jp }
 
\begin{document}
 \maketitle

\clearpage

\vspace{5mm}



\vspace{5mm}

The contact author for correspondence: 

Iwao Sato 

Oyama National College of Technology, 
Oyama, Tochigi 323-0806, JAPAN

E-mail: isato@oyama-ct.ac.jp

\clearpage

\begin{abstract} 
We introduce a generalized Grover matrix of a graph and present an explicit formula for its characteristic polynomial. 
As a corollary, we give the spectra for the generalized Grover matrix of a regular graph. 
Next, we define a zeta function and a generalized zeta function of a graph $G$ with respect to its generalized Grover matrix 
as an analog of the Ihara zeta function and present explicit formulas for their zeta functions for a vertex-transitive graph. 
As applications, we express the limit on the generalized zeta functions of a family of finite vertex-transitive regular graphs 
by an integral. 
Furthermore, we give the limit on the generalized zeta functions of a family of finite tori as an integral expression. 
\end{abstract}

\section{Introduction}

Starting from $p$-adic Selberg zeta functions, Ihara \cite{Ihara} introduced the Ihara zeta functions 
of graphs. 
Bass \cite{Bass} generalized Ihara's result on the Ihara zeta function of 
a regular graph to an irregular graph and showed that its reciprocal is 
a polynomial. 

The Ihara zeta function of a finite graph was extended to an infinite graph. 
Clair \cite{Clair} computed the Ihara zeta function for the infinite grid by using elliptic integrals 
and theta functions. 
Chinta, Jorgenson and Karlsson \cite{Ch} gave a generalized version of 
the determinant formula for the Ihara zeta function associated to 
finite or infinite regular graphs.

There are exciting developments between quantum walk on a graph \cite{Ambainis2003, Kempe2003, 
Kendon2007, Konno2008b, Po, VA} and the Ihara zeta function of a graph. 
 
Ren et al. \cite{Ren} gave a relationship between the discrete-time quantum walk and the Ihara zeta function of a graph.  
Konno and Sato \cite{KS} obtained a formula of the characteristic polynomial of the Grover matrix 
by using the determinant expression for the second weighted zeta function of a graph. 
Komatsu, Konno and Sato \cite{KKS} treated the generalized Ihara zeta function of $ \mathbb{Z} $ as a limit of the Ihara zeta function 
of the cycle graph $C_n $ with $n$ vertices. 

In Grover/Zeta Correspondence \cite{K1}, Komatsu, Konno and Sato defined a zeta function and a generalized zeta function of a graph $G$ 
with respect to its Grover matrix, and presented the limits on the generalized zeta functions and the generalized Ihara zeta functions 
of a family of finite regular graphs as an integral expression by using the Konno-Sato theorem \cite{KS}. 
This result contained the result on the generalized Ihara zeta function in Chinta et al. \cite{Ch}. 
Furthermore, they obtained the limit on the generalized Ihara zeta functions of a family of finite torus as an integral expression, 
and this result contained the result on the Ihara zeta function of the two-dimensional integer lattice $\mathbb{Z}^2$ in Clair \cite{Clair}. 

In Walk/Zeta Correspondence \cite{K2}, Komatsu, Konno and Sato defined a walk-type zeta function without use of the determinant expressions 
of the zeta function of a graph $G$, and presented various properties of walk-type zeta functions of random walk (RW), correlated random walk (CRW) 
quantum walk (QW) and open quantum walk (OQW) on $G$. 
Also, their limit formulas by using integral expressions were presented. 
Konno and Tamura \cite{KT} computed the walk-type zeta functions for the three- and four-state quantum walk, correlated random walk, the multi-state 
random walk on the one-dimensional torus, and the four-state quantum walk, correlated random walk on the two-dimensional torus. 
Furthermore, they gave an extension of the Konno-Sato theorem.

In this paper, we define a generalized Grover matrix of a graph and treat a walk-type zeta function of a vertex-transitive graph with respect to 
its generalized Grover matrix. 

In Section 2, we review the Ihara zeta function of a finite graph and the generalized Ihara zeta function 
of a finite or infinite vertex-transitive graph. 
In Section 3, we state about the time evolution matrix, i.e., the Grover matrix of the Grover walk on a graph. 
In Section 4, we introduce a generalized Grover matrix of a graph and present an explicit formula for its characteristic polynomial. 
As a corollary, we give the spectra for the generalized Grover matrix of a regular graph. 
In Section 5, we define a zeta function and a generalized zeta function of a graph $G$ with respect to its generalized Grover matrix 
as an analog of the Ihara zeta function and present explicit formulas for their zeta functions of a vertex-transitive graph. 
In Section 6, we express the limit on the generalized zeta functions of a family of finite vertex-transitive regular graphs 
by an integral. 
In Section 7, we give the limit on the generalized functions of a family of finite tori as an integral expression. 

\section{The Ihara zeta function of a graph}

All graphs in this paper are assumed to be simple. 
Let $G=(V(G),E(G))$ be a connected graph (without multiple edges and loops) 
with the set $V(G)$ of vertices and the set $E(G)$ of unoriented edges $uv$ 
joining two vertices $u$ and $v$.
Furthermore, let $n=|V(G)|$ and $m=|E(G)|$ be the number of vertices and edges of $G$, respectively. 
For $uv \in E(G)$, an arc $(u,v)$ is the oriented edge (or the arc) from $u$ to $v$. 
Let $D_G$ the symmetric digraph corresponding to $G$. 
Set $D(G)= \{ (u,v),(v,u) \mid uv \in E(G) \} $. 
For $e=(u,v) \in D(G)$, set $u=o(e)$ and $v=t(e)$. 
Furthermore, let $e^{-1}=(v,u)$ be the {\em inverse} of $e=(u,v)$. 
For $v \in V(G)$, the {\em degree} $\deg {}_G v = \deg v = d_v $ of $v$ is the number of vertices 
adjacent to $v$ in $G$.  

A {\em path $P$ of length $n$} in $G$ is a sequence 
$P=(e_1, \ldots ,e_n )$ of $n$ arcs such that $e_i \in D(G)$,
$t( e_i )=o( e_{i+1} ) (1 \leq i \leq n-1)$. 
If $e_i =( v_{i-1} , v_i )$ for $i=1, \ldots , n$, then we write 
$P=(v_0, v_1, \ldots ,v_{n-1}, v_n )$. 
Set $ \mid P \mid =n$, $o(P)=o( e_1 )$ and $t(P)=t( e_n )$. 
Also, $P$ is called an {\em $(o(P),t(P))$-path}. 
We say that a path $P=( e_1 , \ldots , e_n )$ has a {\em backtracking} 
if $ e^{-1}_{i+1} =e_i $ for some $i(1 \leq i \leq n-1)$. 
A $(v, w)$-path is called a {\em $v$-cycle} 
(or {\em $v$-closed path}) if $v=w$. 
The {\em inverse cycle} of a cycle 
$C=( e_1, \ldots ,e_n )$ is the cycle 
$C^{-1} =( e^{-1}_n, \ldots ,e^{-1}_1 )$.

We introduce an equivalence relation between cycles. 
Two cycles $C_1 =(e_1, \ldots ,e_m )$ and 
$C_2 =(f_1, \ldots ,f_m )$ are called {\em equivalent} if 
$f_j =e_{j+k} $ for all $j$. 
The inverse cycle of $C$ is in general not equivalent to $C$. 
Let $[C]$ be the equivalence class which contains a cycle $C$. 
Let $B^r$ be the cycle obtained by going $r$ times around a cycle $B$. 
Such a cycle is called a {\em multiple} of $B$. 
A cycle $C$ is {\em reduced} if 
both $C$ and $C^2 $ have no backtracking. 
Furthermore, a cycle $C$ is {\em prime} if it is not a multiple of 
a strictly smaller cycle. 
Note that each equivalence class of prime, reduced cycles of a graph $G$ 
corresponds to a unique conjugacy class of 
the fundamental group $ \pi {}_1 (G,v)$ of $G$ at a vertex $v$ of $G$. 

The {\em Ihara zeta function} of a graph $G$ is 
a function of a complex variable $u$ with $|u|$ 
sufficiently small, defined by 
\[
{\bf Z} (G, u)= {\bf Z}_G (u)= \prod_{[C]} 
(1- u^{ \mid C \mid } )^{-1} ,
\]
where $[C]$ runs over all equivalence classes of prime, reduced cycles of $G$.

Let $G$ be a connected graph with $n$ vertices $v_1, \ldots ,v_n $. 
The {\em adjacency matrix} ${\bf A}= {\bf A} (G)=(a_{ij} )$ is 
the square matrix such that $a_{ij} =1$ if $v_i$ and $v_j$ are adjacent, 
and $a_{ij} =0$ otherwise.  
If $ \deg {}_G v=k$ (constant) for each $v \in V(G)$, then $G$ is called 
{\em $k$-regular}.

\newtheorem{theorem}{Theorem}
\begin{theorem}[Ihara; Bass] 
Let $G$ be a connected graph. 
Then the reciprocal of the Ihara zeta function of $G$ is given by 
\[
{\bf Z} (G,u )^{-1} =(1- u^2 )^{r-1} 
\det ( {\bf I} -u {\bf A} (G)+ u^2 ( {\bf D} - {\bf I} )) , 
\]
where $r$ is the Betti number of $G$, 
and ${\bf D} =( d_{ij} )$ is the diagonal matrix 
with $d_{ii} = \deg v_i$ and $d_{ij} =0, i \neq j , 
(V(G)= \{ v_1 , \ldots , v_n \} )$. 
\end{theorem}

Let $G=(V(G),E(G))$ be a connected graph and $ x_0 \in V(G)$ 
a fixed vertex.  
Then the {\em generalized Ihara zeta function} $\zeta {}_G (u)$ of $G$ is defined by 
\[
\zeta {}_G (u)= \exp \left( \sum^{\infty}_{m=1} \frac{N^0_m }{m} u^m \right) , 
\]
where $N^0_m $ is the number of reduced $x_0$-cycles of length $m$ in $G$. 
A graph $G$ is called {\em vertex-transitive} if there exists an automorphism $ \phi $ of the automorphism group 
${\rm Aut} (G)$ of $G$ such that $ \phi (u)=v$ for each $u,v \in V(G)$. 
Note that, for a finite vertex-transitive graph, the classical Ihara zeta function is just the above
Ihara zeta function raised to the power equaling the number of vertices: 
\[
\zeta {}_G (u)= {\bf Z} (G,u)^{1/n} . 
\]
Furthermore, the {\em Laplacian} of $G$ is given by 
\[
\Delta = \Delta (G) = {\bf D} - {\bf A} (G). 
\]

A formula for the generalized Ihara zeta function of a vertex-transitive graph is given as follows:

\begin{theorem}[Chinta, Jorgenson and Karlsson] 
Let $G$ be a  vertex-transitive $(q+1)$-regular graph with spectral measure $\mu {}_{\Delta }$ for the Laplacian $\Delta $. 
Then 
\[
\zeta {}_G (u)^{-1} =(1-u^2 )^{(q-1)/2} \exp \left( \int \log (1-(q+1- \lambda )u+q u^2 ) d \mu {}_{\Delta } ( \lambda ) \right) . 
\]
\end{theorem}

Note, if $G$ is a vertex-transitive graph with $n$ vertices, then

\section{The Grover matrix of a graph} 

We define the Grover matrix which is the time evolution matrix of the Grover walk on a graph. 

Let $G$ be a connected graph with $n$ vertices and $m$ edges. 
Then the {\em Grover matrix} ${\bf U} ={\bf U} (G)=( U_{ef} )_{e,f \in D(G)} $ 
of $G$ is defined by 
\[
U_{ef} =\left\{
\begin{array}{ll}
2/d_{t(f)} (=2/d_{o(e)} ) & \mbox{if $t(f)=o(e)$ and $f \neq e^{-1} $, } \\
2/d_{t(f)} -1 & \mbox{if $f= e^{-1} $, } \\
0 & \mbox{otherwise. }
\end{array}
\right. 
\]
The discrete-time quantum walk with the matrix ${\bf U} $ as a time evolution matrix 
is called the {\em Grover walk} on $G$. 
Furthermore, we introduce the {\em positive support} ${\bf F}^+ =( F^+_{ij} )$ of 
a real square matrix ${\bf F} =( F_{ij} )$ as follows: 
\[
F^+_{ij} =\left\{
\begin{array}{ll}
1 & \mbox{if $F_{ij} >0$, } \\
0 & \mbox{otherwise. }
\end{array}
\right.
\]

In Konno and Sato \cite{KS}, they presented the following result. 
The $n \times n$ matrix ${\bf P} = {\bf P} (G)=( P_{uv} )_{u,v \in V(G)}$ is given as follows: 
\[
P_{uv} =\left\{
\begin{array}{ll}
1/( \deg {}_G \ u)  & \mbox{if $(u,v) \in D(G)$, } \\
0 & \mbox{otherwise.}
\end{array}
\right.
\] 
Note that the matrix ${\bf P} (G)$ is the transition probability matrix of the simple random walk on $G$.

\begin{theorem}[Konno and Sato]
Let $G$ be a connected graph with $n$ vertices $v_1 ,\ldots , v_n $ and $m$ edges. 
Then  

\begin{align*}
\det ( {\bf I}_{2m} -u {\bf U} ) &= (1-u^2)^{m-n} \det ((1+u^2) {\bf I}_{n} -2u {\bf P} (G))  
\\
&= \frac{(1-u^2)^{m-n} }{ \prod^n_{i=1} \deg v_i } 
\det ((1+u^2) {\bf D} -2u {\bf A} (G)) . 
\end{align*}
\end{theorem}

This theorem is called the {\em Konno-Sato theorem} (see \cite{KI, Morita}, for example).

Konno and Tamura \cite{KT} extended the Grover matrix. 
Let $G$ be a connected graph with $m$ edges, and $a \in [0,1]$.  
Then the extension ${\bf U}_a ={\bf U}_a (G)=( U^{(a)}_{ef} )_{e,f \in D(G)} $ of the Grover matrix of $G$ 
is defined as follows:  
\[
U^{(a)}_{ef} =\left\{
\begin{array}{ll}
(2/d_{t(f)} -1)a+1 & \mbox{if $t(f)=o(e)$ and $f \neq e^{-1} $, } \\
(2/d_{t(f)} -1)a  & \mbox{if $f= e^{-1} $, } \\
0 & \mbox{otherwise. }
\end{array}
\right. 
\]
If $a=1$, then ${\bf U}_1 = {\bf U}$ is the Grover matrix of $G$. 
In the case of $a=0$, ${\bf U}_0 = {\bf U}^+$ is the positive support of the Grover matrix of $G$. 
Thus, the matrix ${\bf U}_a $ is an extension of the Grover matrix ${\bf U}$ of $G$.

Konno and Tamura \cite{KT} presented the following result for the extension ${\bf U}_a $ of the Grover matrix of a regular graph.

\begin{theorem}[Konno and Tamura]
Let $G$ be a connected $(q+1)$-regular graph with $n$ vertices and $m$ edges. 
Then  
\[
\det ( {\bf I}_{2m} -u {\bf U}_a )= (1-u^2)^{m-n} \det (\{ 1+(q+(1-q)a) u^2 \} {\bf I}_{n} - \{1+q+(1-q)u \} {\bf P} (G)) . 
\] 
\end{theorem}

\section{A generalized Grover matrix of a graph} 

We introduce a generalized Grover matrix of a graph. 

Let $G$ be a connected graph with $m$ edges, $a \in [0,1]$ and $b \in \mathbb{R}$.   
Then a {\em generalized Grover matrix} $\tilde{{\bf U}} = \tilde{{\bf U}} (G)=( \tilde{U}_{ef} )_{e,f \in D(G)} $ of $G$ 
is defined as follows:  
\[
\tilde{U}_{ef} =\left\{
\begin{array}{ll}
(2/d_{t(f)} -1)a+b & \mbox{{\rm if} $t(f)=o(e)$ and $f \neq e^{-1} $, } \\
(2/d_{t(f)} -1)a  & \mbox{if $f= e^{-1} $, } \\
0 & \mbox{{\rm otherwise}. }
\end{array}
\right. 
\]
If $a=b=1$, then $\tilde{{\bf U}} = {\bf U}$ is the Grover matrix of $G$. 
In the case of $a=0$ and $b=1$, $\tilde{{\bf U}} = {\bf U}^+ $ is the positive support of the Grover matrix of $G$. 
Thus, the generalized Grover matrix $\tilde{{\bf U}} $ is a generalization of the Grover matrix ${\bf U}$ and 
the positive support of the Grover matrix ${\bf U}^+$ of $G$.

We present a generalization of the Konno-Sato theorem as follows.

\begin{theorem}[A generalization of the Konno-Sato theorem]
Let $G$ be a connected graph with $n$ vertices $v_1 ,\ldots , v_n $ and $m$ edges, $a \in [0,1]$ and $b \in \mathbb{R}$. 
Then 
\[ 
\det ( {\bf I}_{2m} -u \tilde{{\bf U}} )= \frac{(1- b^2 u^2)^{m-n} }{ \prod^n_{i=1} \deg v_i } 
\det ( {\bf D} \{ (1+b(2a-b) u^2) {\bf I}_n +b(b-a) u^2 {\bf D} \} -u {\bf A}_d ) ,  
\] 
where ${\bf A}_d =( A^{(d)}_{uv} )_{u,v \in V(G)} $ is given as follows: 
\[
A^{(d)}_{uv} =\left\{
\begin{array}{ll}
(2- \deg u)a+b \deg u & \mbox{{\rm if} $(u,v) \in D(G)$, } \\
0 & \mbox{{\rm otherwise}. }
\end{array}
\right. 
\]
\end{theorem}

{\bf Proof }.  At first, let 
\[
w(u,v)= \left( \frac{2}{ \deg u} -1 \right) a+b \ {\rm if} \ (u,v) \in D(G) ,  
\]
and two $n \times n$ matrices ${\bf W} =( w_{uv} )_{u,v \in V(G)} $ and 
${\bf D}_w =( d^{(w)}_{uv} )_{u,v \in V(G)} $ be defined as follows: 
\[
w_{uv} = \left\{
\begin{array}{ll}
w(u,v) & \mbox{if $(u,v) \in D(G)$, } \\
0 & \mbox{otherwise, }
\end{array}
\right.
\
d^{(w)}_{uv} = \left\{
\begin{array}{ll}
\sum_{o(e)=u} w(e) & \mbox{if $u=v$, } \\
0 & \mbox{otherwise. }
\end{array}
\right.
\] 
Furthermore, we define two $2m \times 2m$ matrices ${\bf B}_w = {\bf B}_w (G) = ( B_{ef} )_{e,f \in D(G)} $ and 
${\bf J}_0 =( J_{ef} )_{e,f \in D(G)} $ as follows: 
\[
B_{ef} = \left\{
\begin{array}{ll}
w(f) & \mbox{if $t(e)=o(f)$, } \\
0 & \mbox{otherwise, }
\end{array}
\right.
\
J_{ef} = \left\{
\begin{array}{ll}
1 & \mbox{if $f= e^{-1} $, } \\
0 & \mbox{otherwise. }
\end{array}
\right.
\] 
Then we have 
\[
{}^t \tilde{{\bf U}} = {\bf B}_w - b {\bf J}_0 . 
\]

Next, we introduce $2m \times n$ matrices ${\bf K} =( K_{ev} )_{e \in D(G); v \in V(G)} $ and 
${\bf L} =( L_{ev} )_{e \in D(G); v \in V(G)} $ as follows: 
\[
K_{ev} = \left\{
\begin{array}{ll}
1 & \mbox{if $t(e)=v$, } \\
0 & \mbox{otherwise,  }
\end{array}
\right. 
\ 
L_{ev} = \left\{
\begin{array}{ll}
w(e) & \mbox{if $o(e)=v$, } \\
0 & \mbox{otherwise. }
\end{array}
\right.
\] 
Furthermore, a $2m \times n$ matrix ${\bf M} =( M_{ev} )_{e \in D(G); v \in V(G)} $ be given as follows: 
\[
M_{ev} = \left\{
\begin{array}{ll}
1 & \mbox{if $o(e)=v$, } \\
0 & \mbox{otherwise. }
\end{array}
\right.
\] 
Then we have 
\[
{\bf M} = {\bf J}_0 {\bf K} , \ {\bf K} = {\bf J}_0 {\bf M} . 
\]
Furthermore, we have 
\[
{\bf K} \ {}^t {\bf L} = {\bf B}_w , \ {}^t {\bf L} {\bf K} ={\bf W} , \ {}^t {\bf L} {\bf M} = {\bf D}_w , 
\ {}^t {\bf M} {\bf M} = {}^t {\bf L} {\bf L} = {\bf D} .  
\]

Now, we have 
\begin{align*}
\det ( {\bf I}_{2m} -u \tilde{{\bf U}} ) 
&= \det ( {\bf I}_{2m} -u \ {}^t \tilde{{\bf U}} ) 
\\
&= \det ( {\bf I}_{2m} -u( {\bf B}_w - b {\bf J}_0 ))
\\
&= \det ( {\bf I}_{2m} -u( {\bf K} \ {}^t {\bf L} - b {\bf J}_0 ))
\\ 
&= \det ( {\bf I}_{2m} +bu {\bf J}_0 -u {\bf K} \ {}^t {\bf L} )
\\
&= \det ( {\bf I}_{2m} -u {\bf K} \ {}^t {\bf L} ( {\bf I}_{2m} +bu {\bf J}_0)^{-1} )  
\det ( {\bf I}_{2m} +bu {\bf J}_0 ) . 
\end{align*}
If ${\bf A}$ and ${\bf B}$ are an $m \times n $ and $n \times m$ 
matrices, respectively, then we have 
\[
\det ( {\bf I}_{m} - {\bf A} {\bf B} )= 
\det ( {\bf I}_n - {\bf B} {\bf A} ) . 
\]
Thus, we have 
\[
\det ( {\bf I}_{2m} -u \tilde{{\bf U}} )= \det ( {\bf I}_{n} -u \ {}^t {\bf L} ( {\bf I}_{2m} +bu {\bf J}_0)^{-1} {\bf K} )  
\det ( {\bf I}_{2m} +bu {\bf J}_0 ) . 
\]

But, we have 
\begin{align*}
\det ({\bf I}_{2m} +bu {\bf J}_0 )
&= \det  
\left[ 
\begin{array}{cc}
{\bf I}_m & bu {\bf I}_m \\ 
bu {\bf I}_m & {\bf I}_m 
\end{array} 
\right]
\cdot  
\det  
\left[ 
\begin{array}{cc}
{\bf I}_m & -bu {\bf I}_m \\ 
{\bf 0}_m & {\bf I}_m 
\end{array} 
\right]
\\
&= \det  
\left[ 
\begin{array}{cc}
{\bf I}_m & {\bf 0}_m \\ 
bu {\bf I}_m & {\bf I}_m - b^2 u^2 {\bf I}_m  
\end{array} 
\right] 
\\
&= (1- b^2 u^2 )^m . 
\end{align*}
Furthermore, 
\begin{align*}
({\bf I}_{2m} +bu {\bf J}_0 )^{-1} 
&= \left[ 
\begin{array}{cc}
{\bf I}_m & bu {\bf I}_m \\ 
bu {\bf I}_m & {\bf I}_m 
\end{array} 
\right]^{-1} 
\\
& \sim 
\left[ 
\begin{array}{cccccc}
1 & bu & \ldots &   & 0 \\
bu & 1 &  &  &  \\
\vdots &  & \ddots &  &  \\ 
  &  &  & 1 & bu \\ 
0 &  &  & bu & 1  
\end{array} 
\right]^{-1} 
\\
&= \frac{1}{1- b^2 u^2}  
\left[ 
\begin{array}{ccccc}
1 & -bu & \ldots &  & 0 \\
-bu & 1 &  &  &  \\
\vdots &  & \ddots &  &  \\ 
  &  &  & 1 & -bu \\ 
0 &  &  & -bu & 1   
\end{array} 
\right] 
\\
&= \frac{1}{1- b^2 u^2} ( {\bf I}_{2m} -bu {\bf J}_0 ) . 
\end{align*}

Therefore, it follows that 
\begin{align*}
\det ( {\bf I}_{2m} -u \tilde{{\bf U}} ) 
&=(1- b^2 u^2 )^{m} 
\det \left( {\bf I}_{n} - \frac{u}{1- b^2 u^2 } \ {}^t {\bf L} ( {\bf I}_{2m} -bu {\bf J}_0) {\bf K} \right) 
\\
&=(1- b^2 u^2 )^{m-n} 
\det ((1- b^2 u^2 ) {\bf I}_{n} -u \ {}^t {\bf L} {\bf K} +b u^2 \ {}^t {\bf L} {\bf J}_0 {\bf K} )
\\
&=(1- b^2 u^2 )^{m-n} 
\det ((1- b^2 u^2 ) {\bf I}_{n} -u {\bf W} +b u^2 \ {}^t {\bf L} {\bf M} )
\\
&=(1- b^2 u^2 )^{m-n} 
\det ((1- b^2 u^2 ) {\bf I}_{n} -u {\bf W} +b u^2 {\bf D}_w ) 
\\
&=(1- b^2 u^2 )^{m-n} 
\det ( {\bf I}_{n} -u {\bf W} +b u^2 ( {\bf D}_w -b {\bf I}_n )) . 
\end{align*}

The entries of two matrices ${\bf W} $ and ${\bf D}_w $ are given as follows: 
\[( {\bf W} )_{uv} = \left( \frac{2}{ \deg u} -1 \right) a+b= \frac{1}{ \deg u} \{ (2- \deg u)a +b \deg u \} 
\ {\rm if} \ (u,v) \in D(G) 
\]
and 
\[ 
( {\bf D} )_{uu} = \left\{ \left( \frac{2}{ \deg u} -1 \right) a+b \right\} \deg u =(b-a) \deg u +2a . 
\]
Thus, we have 
\[
{\bf W} = {\bf D}^{-1} {\bf A}_d \ and \ {\bf D}_w =(b-a) {\bf D} +2a {\bf I}_n . 
\]
Therefore, it follows that 
\begin{align*}
\det ( {\bf I}_{2m} -u \tilde{{\bf U}} ) 
&=(1- b^2 u^2 )^{m-n} 
\det ( {\bf I}_{n} -u {\bf D}^{-1} {\bf A}_d +b u^2 ((b-a) {\bf D} +(2a-b) {\bf I}_n ))
\\
&=(1- b^2 u^2 )^{m-n} \det ( {\bf D}^{-1} )  
\det ( {\bf D} -u {\bf A}_d +b u^2 ((b-a) {\bf D}^2 +(2a-b) {\bf D} ))
\\
&= \frac{(1- b^2 u^2)^{m-n} }{ \prod^n_{i=1} \deg v_i } 
\det ( {\bf D} \{ (1+b(2a-b) u^2) {\bf I}_n +b(b-a) u^2 {\bf D} \} -u {\bf A}_d ) .   
\end{align*}
$\Box$

For a $(q+1)$-regular graph, we obtain the following result.

\newtheorem{corollary}{Corollary}
\begin{corollary}
Let $G$ be a connected $(q+1)$-regular graph with $n$ vertices $v_1 ,\ldots , v_n $ and $m$ edges, $a \in [0,1]$ and $b \in \mathbb{R}$. 
Then 
\begin{align*}
& \det ( {\bf I}_{2m} -u \tilde{{\bf U}} )
\\
&=(1- b^2 u^2)^{m-n}   
\det ( \{ 1+b((1-q)a+bq) u^2 \} {\bf I}_n -u ((1-q)a+b(q+1)) {\bf P} (G)) . 
\end{align*}
\end{corollary}

{\bf Proof }.  At first, we have 
\[
{\bf D} =(q+1) {\bf I}_n , \ {\bf P} = \frac{1}{q+1} {\bf A} , \ {\bf A}_d = \{ (1-q)a+b(q+1) \} {\bf A} . 
\]
Thus, 
\begin{align*}
& \det ( {\bf I}_{2m} -u \tilde{{\bf U}} ) 
\\
&= \frac{(1- b^2 u^2)^{m-n} }{ (q+1)^n } 
\det ((q+1)(1+b(2a-b) u^2 )+(q+1)(b-a)b u^2 ) {\bf I}_n -u((1-q)a+b(q+1)) {\bf A} ) 
\\
&=(1- b^2 u^2)^{m-n}   
\det ( \{ 1+b((1-q)a+bq) u^2 \} {\bf I}_n -u ((1-q)a+b(q+1)) {\bf P} (G)) . 
\end{align*}
$\Box$

If $G$ is a $(q+1)$-regular graph with $n$ vertices, then we have 
\[
{\bf P} = \frac{1}{q+1} {\bf A} = \frac{1}{q+1} ( {\bf D} - \Delta )= {\bf I}_n - \frac{1}{q+1} \Delta . 
\]
By Corollary 1, we obtain the following result.

\begin{corollary}
Let $G$ be a connected $(q+1)$-regular graph with $n$ vertices $v_1 ,\ldots , v_n $ and $m$ edges, $a \in [0,1]$ and $b \in \mathbb{R}$. 
Then   
\begin{align*}
& \det ( {\bf I}_{2m} -u \tilde{{\bf U}} ) 
\\
&=(1- b^2 u^2)^{m-n}   
\det \bigg( \{ 1-u((1-q)a+b(q+1))+b((1-q)a+bq) u^2 \} {\bf I}_n 
\\ 
& +u \bigg( b- \frac{q-1}{q+1} a \bigg) \Delta \bigg) . 
\end{align*}
\end{corollary}

{\bf Proof}.  By Corollary 1, we have 
\begin{align*}
& \det ( {\bf I}_{2m} -u \tilde{{\bf U}} )
\\
&=(1- b^2 u^2)^{m-n}   
\det ( \{ 1+b((1-q)a+bq) u^2 \} {\bf I}_n 
\\ 
& -u ((1-q)a+b(q+1))( {\bf I}_n - \frac{1}{q+1} \Delta ))
\\
&=(1- b^2 u^2)^{m-n}   
\det \bigg( \{ 1-u((1-q)a+b(q+1))+b((1-q)a+bq) u^2 \} {\bf I}_n 
\\
& +u \bigg( b- \frac{q-1}{q+1} a \bigg) \Delta \bigg) .
\end{align*}
$\Box$

Substituting $u=1 /\lambda $ to Corollaries 1 and 2, we obtain the following result.

\begin{corollary}
Let $G$ be a connected $(q+1)$-regular graph with $n$ vertices $v_1 ,\ldots , v_n $ and $m$ edges, $a \in [0,1]$ and $b \in \mathbb{R}$. 
Then  
\begin{align*} 
& \det ( \lambda {\bf I}_{2m} - \tilde{{\bf U}} ) 
\\
&=( \lambda - b^2 )^{m-n}   
\det ( \{ \lambda {}^2 +b((1-q)a+bq) \} {\bf I}_n - \lambda ((1-q)a+b(q+1)) {\bf P} (G)) 
\\
&= ( \lambda - b^2 )^{m-n}   
\det \bigg( \{ \lambda {}^2 - \lambda ((1-q)a+b(q+1))+b((1-q)a+bq) \} {\bf I}_n 
\\
& + \lambda \bigg( b- \frac{q-1}{q+1} a \bigg) \Delta \bigg) .
\end{align*}
\end{corollary}

By Corollary 3, the following result holds. 
Let ${\rm Spec} ( {\bf F} )$ be the set of eigenvalues of a square matrix ${\bf F} $.

\begin{corollary}
Let $G$ be a connected $(q+1)$-regular graph with $n$ vertices $v_1 ,\ldots , v_n $ and $m$ edges, $a \in [0,1]$ 
and $b \in \mathbb{R}$. 
Set $\eta =(1-q)a+b(q+1)$ and $ \sigma =b((1-q)a+bq)$.  
Then the spectra of the generalized Grover matrix $\tilde{{\bf U}} $ are given as follows: 
\begin{enumerate} 
\item $2n$ eigenvalues: 
\[
\lambda = \frac{ \mu \eta \pm \sqrt{ \mu {}^2 \eta {}^2 -4 \sigma } }{2} , \ \mu \in {\rm Spec} ( {\bf P} ); 
\]
\item $2(m-n)$ eigenvalues: $\pm b$ with multiplicities $m-n$. 
\end{enumerate} 
\end{corollary}

{\bf Proof }. By Corollary 3, we have 
\[ 
\det ( \lambda {\bf I}_{2m} - \tilde{{\bf U}} )=( \lambda - b^2 )^{m-n}   
\prod_{ \mu \in {\rm Spec} ( {\bf P} )} ( \lambda {}^2 + \sigma - \mu \eta \lambda ) . 
\]
Solving $ \lambda {}^2 -\mu \eta \lambda + \sigma =0$, we obtain 
\[
\lambda =\frac{ \mu \eta \pm \sqrt{ \mu {}^2 \eta {}^2 -4 \sigma }}{2} . 
\] 
The result follows. 
$\Box$

\section{A generalized Grover/Zeta Correspondence} 

Now, we propose a new zeta function of a graph.  
Let $G$ be a connected graph with $m$ edges, $a \in [0,1]$ and $b \in \mathbb{R}$.  
Then we define the {\em $(a,b)$-zeta function} $ {\bf Z}_{a,b} (G, u)$ of $G$ is defined as follows:    
\[
{\bf Z}_{a,b} (G, u)^{-1} = {\bf Z}_{a,b} (u)^{-1} = \det ( {\bf I}_{2m} -u \tilde{{\bf U}} ) .    
\]

By Corollaries 1 and 2, we obtain the following result.

\newtheorem{proposition}{Proposition}
\begin{proposition}
Let $G$ be a connected $(q+1)$-regular graph with $n$ vertices and $m$ edges, $a \in [0,1]$ and $b \in \mathbb{R}$. 
Set $\eta =(1-q)a+b(q+1)$ and $ \sigma =b((1-q)a+bq)$.  
Then 
\begin{align*}
{\bf Z}_{a,b} (G,u)^{-1} 
&=(1- b^2 u^2)^{m-n}   
\det ( \{ 1+ \sigma u^2 \} {\bf I}_n - \eta u {\bf P} (G)) 
\\
&=(1- b^2 u^2)^{m-n}   
\det \left( \{ 1- \eta u+ \sigma u^2 \} {\bf I}_n + \frac{\eta u}{q+1} \Delta \right) .  
\end{align*}
\end{proposition}

By Theorem 3, we obtain the exponential expression 
for ${\bf Z}_{a,b} (u)$. 
We give a weight functions $w: D(G) \times D(G) \longrightarrow \mathbb{C} $ as follows: 
\[
w (f,e) =\left\{
\begin{array}{ll}
(2/ \deg t(f)-1)a+b  & \mbox{if $t(f)=o(e)$ and $f \neq e^{-1} $, } \\
(2/ \deg t(f) -1)a & \mbox{if $f= e^{-1} $, } \\
0 & \mbox{otherwise. }
\end{array}
\right.
\] 
For a cycle $C=( e_1, e_2 , \ldots , e_r )$, let 
\[
w(C)=w(e_1 , e_2 ) \cdots w( e_{r-1} , e_r) w( e_r , e_1 ) . 
\]

\begin{theorem} 
Let $G$ be a connected graph with $m$ edges, $a \in [0,1]$ and $b \in \mathbb{R}$.  
Then 
\[
{\bf Z}_{a,b} (u)= \exp \left(\sum^{\infty}_{r=1} \frac{N_r}{r} u^r \right) ,  
\]
where $N_r $ is defined by 
\[
N_r = \sum \{ w(C) \mid C: \ a \ cycle \ of \ length \ r \ in \ G \} . 
\] 
\end{theorem}

{\bf Proof}. By the definition of ${\bf Z}_{a,b} (u)$, we have 
\[
\log {\bf Z}_{a,b} (u)=\log \{ \det ( {\bf I}_{2m} -u \tilde{{\bf U}} )^{-1} \} 
=- {\rm Tr} [ \log ( {\bf I}_{2m} -u \tilde{{\bf U}} )] 
= \sum^{\infty}_{r=1} \frac{1}{r} {\rm Tr} [ \tilde{{\bf U}}^r ] u^r . 
\]
Since  
\[
w(f,e)=( \tilde{{\bf U}} )_{ef} , \ e,f \in D(G),  
\]
we have 
\[
{\rm Tr} [ \tilde{{\bf U}}^r ]= \sum \{ w(C) \mid C: \ a \ cycle \ of \ length \ r \ in \ G \} = N_r . 
\]
Hence, 
\[
\log {\bf Z}_{a,b} (u)= \sum^{\infty}_{r=1} \frac{N_r }{r} u^r . 
\]
Thus, 
\[
{\bf Z}_{a,b} (u)= \exp \left(\sum^{\infty}_{r=1} \frac{N_r }{r} u^r \right) . 
\] 
$\Box$

Next, we define a generalized zeta function with respect to the generalized Grover matrix of a graph. 
Let $G=(V(G),E(G))$ be a connected graph, $ x_0 \in V(G)$ a fixed vertex, $a \in [0,1]$ and $b \in \mathbb{R}$.  
Then the {\em generalized $(a,b)$-zeta function} $\zeta_{a,b} (G, u)$ of $G$ is defined by 
\[
\zeta_{a,b} (G, u)= \exp \left( \sum^{\infty}_{r=1} \frac{N^0_r }{r} u^r \right) , 
\]
where 
\[
N^0_r = \sum \{ w(C) \mid C: \ an \ x_0-cycle \ of \ length \ r \ in \ G \} . 
\] 
Note, if $G$ is a vertex-transitive graph with $n$ vertices, then 
\begin{equation} 
\zeta_{a,b} (u)= {\bf Z}_{a,b} (G,u)^{1/n} . 
\end{equation}

If $a=b=1$, then $\zeta_{1,1} (G, u)= \overline{\zeta} (G, u)$ is the generalized zeta function of $G$ (see \cite{K1}). 
In the case of $a=0$ and $b=1$, $\zeta_{0,1} (G, u)= \zeta (G, u)$ is the generalized Ihara zeta function of $G$. 
Thus, the generalized $(a,b)$-zeta function $\zeta_{a,b} (G, u)$ is a generalization of the generalized zeta function 
and the generalized Ihara zeta function of $G$.

Now, we present an explicit formula for the generalized $(a,b)$-zeta function for a vertex-transitive graph. 

Let $G$ be a vertex-transitive $(q+1)$-regular graph with $n$ vertices and $m$ edges. 
Then, since $m=(q+1)n/2$, we have 
\[
\frac{m-n}{n}= \frac{q-1}{2} . 
\]

By Proposition 1, we obtain the following result.

\begin{theorem}[Generalized Grover/Zeta Correspondence] 
Let $G$ be a connected vertex-transitive $(q+1)$-regular graph with $n$ vertices and $m$ edges, 
$a \in [0,1]$ and $b \in \mathbb{R}$.  
Set $\eta =(1-q)a+b(q+1)$ and $ \sigma =b((1-q)a+bq)$.   
Then 
\begin{equation} 
\zeta_{a,b} (G,u)^{-1} =(1- b^2 u^2)^{(q-1)/2} 
\exp \left[ \frac{1}{n} \sum_{ \lambda \in {\rm Spec}( {\bf P} )} \log \{ (1+ \sigma u^2 )- \eta u \lambda \} \right] ,
\end{equation} 
\begin{equation}  
\zeta_{a,b} (G,u)^{-1} =(1- b^2 u^2)^{(q-1)/2}  
\exp \left[ \frac{1}{n} \sum_{ \lambda \in {\rm Spec}( \Delta )} \log \left\{ (1- \eta u+ \sigma u^2 )
+ \frac{\eta u}{q+1} \lambda \right\} \right] . 
\end{equation} 
\end{theorem}

{\bf Proof}.  By (1) and Corollary 1, we have 
\begin{align*}
{\zeta}_{a,b} (G, u)^{-1} 
&= {\bf Z}_{a,b} (G, u)^{-1/n} = \det ( {\bf I}_{2m} -u \tilde{{\bf U}} )^{1/n} 
\\
&=(1- b^2 u^2 )^{(m-n)/n} \{ \det ((1+ \sigma u^2 ) {\bf I}_{n} - \eta u {\bf P} ) \} {}^{1/n} 
\\
&=(1- b^2 u^2 )^{(q-1)/2} \left\{ \prod_{ \lambda \in {\rm Spec}( {\bf P} )} ((1+ \sigma u^2 )- \eta u \lambda ) \right\}^{1/n} 
\\
&=(1- b^2 u^2 )^{(q-1)/2} \exp \left[ \log \left\{ \prod_{ \lambda \in {\rm Spec}( {\bf P} )} 
((1+ \sigma u^2 )- \eta u \lambda )^{1/n} \right\} \right] 
\\
&=(1- b^2 u^2 )^{(q-1)/2} \exp \left[ \frac{1}{n} \sum_{ \lambda \in {\rm Spec}( {\bf P} )} \log \{ (1+ \sigma u^2 )- \eta u \lambda \} \right] . 
\end{align*}

Similarly, the second formula follows. 
$\Box$

\section{The generalized $(a,b)$-zeta functions for the series of regular graphs} 

We present an explicit formula for the generalized $(a,b)$-functions for the series of regular graphs. 
Let $\{ G_n \}^{\infty}_{n=1} $ be a series of finite vertex-transitive $(q+1)$-regular graphs such that 
\[
\lim {}_{n \rightarrow \infty} |V(G_n )|= \infty . 
\]
Then we have 
\[
\frac{|E(G_n )| -|V(G_n )|}{|V(G_n )|}= \frac{(q-1)|V(G_n )|}{2|V(G_n )|}= \frac{q-1}{2} . 
\]
Set 
\[
\nu_n = |V(G_n)|, \ m_n =|E(G_n )| .
\]

Then the following result holds.

\begin{theorem} 
Let $\{ G_n \}^{\infty}_{n=1} $ be a series of finite vertex-transitive $(q+1)$-regular graphs such that 
\[
\lim {}_{n \rightarrow \infty} |V(G_n )|= \infty . 
\]
Furthermore let $a \in [0,1]$, $b \in \mathbb{R}$, $\eta =(1-q)a+b(q+1)$ and $ \sigma =b((1-q)a+bq)$.   
Then
\begin{enumerate}   
\item $\lim_{n \rightarrow \infty} \zeta_{a,b} (G_n , u)^{-1} =
(1-u^2 )^{(q-1)/2} \exp \left[ \int \log \{ (1+ \sigma u^2 )- \eta u \lambda \} d \mu_{P} ( \lambda ) \right] $;  
\item $\lim_{n \rightarrow \infty} \zeta_{a,b} (G_n , u)^{-1} =
(1-u^2 )^{(q-1)/2} \exp \left[ \int \log \{ (1- \eta u+ \sigma u^2 )+ \frac{ \eta u}{q+1} \lambda \} d \mu_{\Delta} ( \lambda ) \right] $, 
\end{enumerate} 
where $d \mu_P ( \lambda )$ and $d \mu_{\Delta } ( \lambda )$ are the spectral measures for the transition operator ${\bf P}$ 
and the Laplacian $\Delta $. 
\end{theorem}

{\bf Proof}.  By Theorem 7, we have 
\[  
\lim_{n \rightarrow \infty} \zeta_{a,b} (G_n , u)^{-1} =
(1-u^2 )^{(q-1)/2} \exp \left[ \int \log \{ (1+ \sigma u^2 )- \eta u \lambda \} d \mu ( \lambda ) \right] .   
\]

Similarly, the second formula follows. 
$\Box$

If $a=b=1$, then we obtain the Grover/Zeta Correspondence (see \cite{K1}).

\begin{corollary}[Grover/Zeta Correspondence]  
Let $\{ G_n \}^{\infty}_{n=1} $ be a series of finite vertex-transitive $(q+1)$-regular graphs such that 
\[
\lim {}_{n \rightarrow \infty} |V(G_n )|= \infty . 
\]
Then
\begin{enumerate} 
\item $\lim_{n \rightarrow \infty} \overline{\zeta}_{G_n} (u)^{-1} =
(1-u^2 )^{(q-1)/2} \exp \left[ \int \log \{ (1+ u^2 )-2u \lambda \} d \mu_P ( \lambda ) \right] $;   
\item $\lim_{n \rightarrow \infty} \overline{\zeta}_{G_n} (u)^{-1} =
(1-u^2 )^{(q-1)/2} \exp \left[ \int \log \{ (1-2u+ u^2 )+ \frac{2u}{q+1} \lambda \} d \mu_{\Delta } ( \lambda ) \right] $,  
\end{enumerate} 
where $d \mu_P ( \lambda )$ and $d \mu_{\Delta } ( \lambda )$ are the spectral measures for the transition operator ${\bf P}$ 
and the Laplacian $\Delta $. 
\end{corollary}

In the case of $a=0$ and $b=1$, we obtain the Grover(Positive Support)/Ihara Zeta Correspondence (see \cite{K1}).

\begin{theorem}[Grover(Positive Support)/Ihara Zeta Correspondence]  
Let $\{ G_n \}^{\infty}_{n=1} $ be a series of finite vertex-transitive $(q+1)$-regular graphs such that 
\[
\lim {}_{n \rightarrow \infty} |V(G_n )|= \infty . 
\]
Then
\begin{enumerate} 
\item $\lim_{n \rightarrow \infty} \zeta_{G_n} (u)^{-1} =
(1-u^2 )^{(q-1)/2} \exp \left[ \int \log \{ (1+q u^2 )-(q+1)u \lambda \} d \mu_P ( \lambda ) \right] $;   
\item $\lim_{n \rightarrow \infty} \zeta_{G_n} (u)^{-1} =
(1-u^2 )^{(q-1)/2} \exp \left[ \int \log \{ (1+q u^2 )-(q+1- \lambda ) u \} d \mu_{\Delta } ( \lambda ) \right] $,  
\end{enumerate} 
where $d \mu_P ( \lambda )$ and $d \mu_{\Delta } ( \lambda )$ are the spectral measures for the transition operator ${\bf P}$ 
and the Laplacian $\Delta $. 
\end{theorem}

The second formula is Theorem 1.3 in Chinta et al. \cite{Ch}.

\section{Torus cases} 

We consider the generalized $(a,b)$-zeta function of the $d$-dimensional integer lattice $\mathbb{Z}^d \ (d \geq 2)$.  

Let $T^d_N \ (d \geq 2)$ be the {\em $d$-dimensional torus} ({\em graph}) with $N^d$ vertices. 
Its vertices are located in coordinates $i_1 , i_2 , \ldots i_d $ of a $d$-dimensional Euclidian space $\mathbb{R}^d $, 
where $i_j \in \{ 0,1, \ldots , N-1 \} $ for any $j$ from 0 to $d-1$. 
A vertex $v$ is adjacent to a vertex $w$ if and only if they have $d-1$ coordinates that are the same, 
and for the remaining coordinate $k$, we have $|i^v_k - i^w_k |=1$, where $i^v_k $ and $i^w_k $ are the $k$-th coordinate 
of $v$ and $w$, respectively. 
Then we have 
\[
|E( T^d_N )|=d N^d ,  
\]
and $T^d_N$ is a vertex-transitive $2d$-regular graph.  

By Corollary 1, we obtain the following result. 
\begin{equation} 
{\bf Z}_{a,b} (T^d_N ,u)^{-1} = \det ( {\bf I}_{2d N^d } -u \tilde{{\bf U}} (T^d_N)) 
=(1- u^2 )^{(d-1) N^d} \det ((1+ \sigma  u^2) {\bf I}_{N^d} - \eta u {\bf P}^{(s)} (T^d_N )) .   
\end{equation} 
Here, it is known that ${\rm Spec}( {\bf P}^{(s)} (T^d_N ))$ is given as follows (see \cite{Spitzer}): 
\[
{\rm Spec}( {\bf P}^{(s)} (T^d_N ))= \left\{ \frac{1}{d} \sum^d_{j=1} \cos \left( \frac{2 \pi k_j }{N} \right) \Bigg|  
k_1 , \ldots , k_d \in \{ 0,1, \ldots , N-1 \} \right\} .
\]
Thus, 
\[  
\displaystyle 
\zeta {}_{a,b} (T^d_N ,u)^{-1} =(1- u^2 )^{d-1} \exp \left[ \frac{1}{N^d } 
\sum^d_{j=1} \sum^{N-1}_{ k_j =0} \log \left\{ (1+ \sigma u^2 )- \frac{ \eta u}{d} \sum^d_{j=1} 
\cos \left( \frac{2 \pi k_j }{N} \right) \right\} \right] . 
\]

Therefore, we obtain the following theorem.

\begin{theorem}[Generalized Grover/Zeta Correspondence ($T^d_N$ case)]    
Let $T^d_N \ (d \geq 2)$ be the $d$-dimensional torus with $N^d$ vertices.  
Furthermore let $a \in [0,1]$, $b \in \mathbb{R}$, $\eta =2(1-d)a+2db$ and $ \sigma =b(2a-b)+2db(b-a)$.   
Then 
\[ 
\displaystyle 
\lim_{n \rightarrow \infty} {\zeta}_{a,b} (T^d_N ,u)^{-1} =(1- u^2 )^{d-1} \exp \left[ \int^{2 \pi}_{0} \dots \int^{2 \pi}_{0}  
\log \left\{ (1+ \sigma u^2 )- \frac{ \eta u}{d} \sum^d_{j=1} \cos \theta_j \right\} 
\frac{d \theta_1}{2 \pi } \cdots \frac{d \theta_d}{2 \pi } \right] ,  
\] 
where $\int^{2 \pi}_{0} \dots \int^{2 \pi}_{0} $ is the $d$-th multiple integral and 
$ \frac{d \theta_1}{2 \pi } \cdots \frac{d \theta_d}{2 \pi } $ is the uniform measure on $[0, 2 \pi )^d $.  
\end{theorem}

If $a=b=1$, then we obtain the Grover/Zeta Correspondence ($T^d_N$ case) (see \cite{K1}).

\begin{corollary}[Grover/Zeta Correspondence ($T^d_N$ case)]      
Let $T^d_N \ (d \geq 2)$ be the $d$-dimensional torus with $N^d$ vertices.  
Then 
\[ 
\displaystyle 
\lim_{n \rightarrow \infty} \overline{{\zeta}} (T^d_N ,u)^{-1} =(1- u^2 )^{d-1} \exp \left[ \int^{2 \pi}_{0} \dots \int^{2 \pi}_{0}  
\log \left\{ (1+ u^2 )- \frac{2u}{d} \sum^d_{j=1} \cos \theta_j \right\} \frac{d \theta_1}{2 \pi } \cdots \frac{d \theta_d}{2 \pi } \right] ,  
\] 
where $\int^{2 \pi}_{0} \dots \int^{2 \pi}_{0} $ is the $d$-th multiple integral and 
$ \frac{d \theta_1}{2 \pi } \cdots \frac{d \theta_d}{2 \pi } $ is the uniform measure on $[0, 2 \pi )^d $.  
\end{corollary}

In the case of $a=0$ and $b=1$, we obtain the Grover(Positive Support)/Ihara Zeta Correspondence ($T^d_N$ case)  (see \cite{K1}).

\begin{corollary}[Grover(Positive Support)/Ihara Zeta Correspondence ($T^d_N$ case )]    
Let $T^d_N \ (d \geq 2)$ be the $d$-dimensional torus with $N^d$ vertices. 
Then  
\[ 
\displaystyle 
\lim_{n \rightarrow \infty} {\zeta} (T^d_N ,u)^{-1} =(1- u^2 )^{d-1} \exp \left[ \int^{2 \pi}_{0} \dots \int^{2 \pi}_{0}  
\log \left\{ (1+(2d-1) u^2 )-2u \sum^d_{j=1} \cos \theta_j \right\} \frac{d \theta_1}{2 \pi } \cdots \frac{d \theta_d}{2 \pi } \right] ,  
\] 
where $\int^{2 \pi}_{0} \dots \int^{2 \pi}_{0} $ is the $d$-th multiple integral and 
$ \frac{d \theta_1}{2 \pi } \cdots \frac{d \theta_d}{2 \pi } $ is the uniform measure on $[0, 2 \pi )^d $.  
\end{corollary}

Specially, in the case of $d=2$, we obtain the following result.

\begin{corollary}     
Let $T^2_N $ be the $2$-dimensional torus with $N^2$ vertices. 
Then  
\[ 
\displaystyle 
\lim_{n \rightarrow \infty} {\zeta} (T^2_N ,u)^{-1} =(1- u^2 ) \exp \left[ \int^{2 \pi}_{0} \int^{2 \pi}_{0}  
\log \left\{ (1+3 u^2 )-2u \sum^2_{j=1} \cos \theta_j \right\} \frac{d \theta_1}{2 \pi } \frac{d \theta_2}{2 \pi } \right] .   
\] 
\end{corollary}

This result corresponds to Equation (10) in Clair \cite{Clair}. 

Finally, we should remark $d=1$ case studied in Komatsu, Konno and Sato \cite{KKS}. 
In this case, we easily check ${\bf U} = {\bf U}^+ $.   
So we can apply both of our results (Corollaries 6 and 7) and get the same result given by Komatsu, Konno and Sato \cite{KKS}.

\end{document}